\documentclass[12pt]{amsart}
\usepackage{enumerate}
\usepackage{amsmath}
\usepackage{amssymb}
\usepackage{amsfonts}
\usepackage{amsthm, upref}
\usepackage{graphicx}
\usepackage[usenames, dvipsnames]{color}

    \numberwithin{equation}{section}

\usepackage{bm}
         
         \bmdefine\alphab{\mathbf{\alpha}}
\bmdefine\betab{\mathbf{\beta}}
\bmdefine\sigmab{\mathbf{\sigma}}

\newcommand{\comment}[1]{}
\newcommand{\eq}{\begin{equation}}
\newcommand{\en}{\end{equation}}
\newcommand{\pp}{\mathbb{P}}
\newcommand{\rr}{\mathbb{R}}

\newcommand{\ev}{\mathbb E}

\newcommand{\ep}{\hfill $\Box$}

\begin{document}

\theoremstyle{plain}
\newtheorem{thm}{Theorem}
\newtheorem{lemma}[thm]{Lemma}
\newtheorem{prop}[thm]{Proposition}
\newtheorem{cor}[thm]{Corollary}

\theoremstyle{definition}
\newtheorem{defn}{Definition}
\newtheorem{cond}{Condition}
\newtheorem{asmp}{Assumption}
\newtheorem{notn}{Notation}
\newtheorem{prb}{Problem}

\theoremstyle{remark}
\newtheorem{rmk}{Remark}
\newtheorem{exm}{Example}
\newtheorem{clm}{Claim}

\title[Transformation between martingales]{On a non-linear transformation between Brownian martingales}
\author{Mykhaylo Shkolnikov}
\address{Mathematical Sciences Research Institute \\ 17 Gauss Way \\ Berkeley, CA 94720-5070}
\email{mshkolni@gmail.com}
\date{\today}

\begin{abstract}
The paper studies a non-linear transformation between Brownian martingales, which is given by the inverse of the pricing operator in the mathematical finance terminology. Subsequently, the solvability of systems of equations corresponding to such transformations is investigated. The latter give rise to novel monotone pathwise couplings of an arbitrary number of certain diffusion processes with varying diffusion coefficients. In the case that there is an uncountable number of these diffusion processes and that the index set is an interval such couplings can be viewed as models for the growth of one-dimensional random surfaces. With this motivation in mind, we derive the appropriate stochastic partial differential equations for the growth of such surfaces.      
\end{abstract}

\maketitle

\section{Introduction}

The starting point of the paper is the equation
\eq\label{cond}
M(t)=\ev[v(X(T))|{\mathcal F}^X(t)],\quad t\in[0,T],
\en
where $T$ is a positive real number, $M(t)$, $t\in[0,T]$ and $X(t)$, $t\in[0,T]$ are stochastic processes with continuous real-valued paths, 
$({\mathcal F}^X(t))_{t\in[0,T]}$ is the filtration generated by $X$ and $v:\,\rr\rightarrow\rr$ is a (deterministic) function. This type of equations is fundamental to the mathematical theory of asset pricing, where the process $X$ stands for the price process of an asset in a financial market, $v(X(T))$ for the payoff of a contract on this asset and $M$ for the price process of this contract. In this context, one usually assumes that the process $X$ is of a particular form (for example, a logarithmic Brownian motion in the classical Black-Scholes model as in \cite{BS}, \cite{Me2}) and tries to determine exactly or to obtain estimates on the process $M$. It is standard to model the asset price process $X$ by a continuous diffusion process and one assumes that $X$ is a martingale, which implies the absence of arbitrage in the market given by $X$ and $M$ by the Fundamental Theorem of Asset Pricing (see \cite{HK}, \cite{HP}, \cite{DR}, \cite{DS} and the references therein). For the rest of the paper, we will call continuous diffusion processes, which are martingales, \textit{Brownian martingales}. 

\bigskip

The starting point of our analysis is the following question asked by David Aldous (\cite{Al}): Suppose that $M$ is a given Brownian martingale and that the function $v$ is given. Can one find a process $X$ in the class of Brownian martingales such that \eqref{cond} holds in the sense of equality of probability laws? In other words, can one find a Brownian martingale $X$ such that the process $\ev[v(X(T))|{\mathcal F}^X(t)]$, $t\in[0,T]$ coincides in law with a given Brownian martingale $M$? The first main result of the paper is that, under certain regularity assumptions on $M$ and $v$, the question can be answered in the affirmative. The operation of solving \eqref{cond} for $X$ can be viewed as the inverse of the pricing operator, which maps $(X,v)$ to $M$ via \eqref{cond}. We remark at this point that, although we allow $M$ and $X$ to take negative values throughout, one can enforce the nonnegativity of $M$ and $X$ suggested by the financial interpretation by letting $M$ be a Brownian martingale taking only nonnegative values and choosing $v$ such that the preimage of the state space of $M$ under $v$ is a subset of $[0,\infty)$. Similarly, one can ensure that $M$ and $X$ take values in the desired bounded subintervals of $[0,\infty)$.     

\bigskip

In the second part of the paper, we study systems of equations
\eq\label{cond_sys}
M^\lambda(\cdot)\overset{d}{=}\ev[v^\lambda(X(T))|{\mathcal F}^X(\cdot)],\quad \lambda\in\Lambda
\en
for arbitrary index sets $\Lambda$. We provide a necessary and sufficient condition on the family $M^\lambda$, $\lambda\in\Lambda$ of Brownian martingales under which it is possible to find a Brownian martingale $X$ and functions $v^\lambda:\,\rr\rightarrow\rr$, $\lambda\in\Lambda$ of certain regularity satisfying \eqref{cond_sys}. Moreover, if $\Lambda$ is endowed with a partial order $\prec$, we give a condition for the system \eqref{cond_sys} to be solvable in the specified sense with functions $v^\lambda:\,\rr\rightarrow\rr$, $\lambda\in\Lambda$, which satisfy $v^{\lambda_1}\leq v^{\lambda_2}$ whenever $\lambda_1\prec\lambda_2$. Clearly, the result are monotone pathwise couplings of the Brownian martingales $M^\lambda$, $\lambda\in\Lambda$. If $\Lambda\subset\rr$ is an interval endowed with the usual total order, then by viewing $\lambda$ as the time variable and $t$ as the space variable, one may regard the process $M^\lambda(t)$, $(\lambda,t)\in\Lambda\times[0,T]$ as a model for the growth of a one-dimensional random surface over time. We show that the growth of the surface can be described by a stochastic partial differential equation (SPDE). 
 
\bigskip

To state our first result rigorously, we introduce some notation. In the context of equation \eqref{cond}, we denote the state space of $M$ by $S\subset\rr$, which we assume to be a (bounded or unbounded) open interval, and write ${\mathcal L}_t=\frac{1}{2}\,a(t,x)\,\frac{\mathrm{d}^2}{\mathrm{d}x^2}$ for the time-dependent generator of the process $M$. Moreover, we let ${\mathcal G}$ be the space of continuously differentiable functions $v:\,\rr\rightarrow\rr$, whose derivatives $v'$ are bounded between two positive constants. 

\begin{thm}\label{thm1}
Suppose that $S$ is contained in the range of $v$ and that the transition densities $p^M(t,x;T,y)$, $t\in[0,T)$, $x,y\in S$ of $M$ exist and possess continuous partial derivatives $p^M_t$, $p^M_x$, $p^M_{xx}$ and $p^M_{xxx}$. Moreover, assume that there are constants $C_1,C_2>0$ such that
\eq\label{pMx_growth}
\int_S |p^M_x(t,x;T,y)|(1+|y|)\,\mathrm{d}y\leq C_1 e^{C_2\,x^2}
\en
for all $(t,x)\in[0,T)\times S$ and that for any given $t\in[0,T)$ and $x\in S$ there exists an open neighborhood $U$ of $(t,x)$ and a function $F:\,S\rightarrow[0,\infty)$ such that
\begin{eqnarray*}
&&\big(|p^M_t|+|p^M_x|+|p^M_{xx}|+|p^M_{xxx}|\big)(\tilde{t},\tilde{x};T,y)\leq F(y),\,(\tilde{t},\tilde{x},y)\in U\times S,\\
&&\int_S F(y)(1+|y|)\,\mathrm{d}y<\infty. 
\end{eqnarray*}
Finally, suppose that $a$ takes only positive values, its partial derivative $a_x$ is continuous and the estimates 
\eq
|a(t,x)|\leq C_3,\quad |a_x(t,x)|\leq C_3(1+|x|),\quad (t,x)\in[0,T]\times S
\en
hold for some constant $C_3>0$. Then, for any $v\in{\mathcal G}$, there exists a unique Brownian martingale $X$ satisfying the equation \eqref{cond} in law and it solves the stochastic differential equation (SDE)  
\eq\label{Xsde}
\mathrm{d}X(t)=h_x(t,h^{(-1)}(t,X(t)))\,\sqrt{a(t,h^{(-1)}(t,X(t)))}\,\mathrm{d}B(t),\quad t\in[0,T],
\en
where $h$ solves the backward Cauchy problem
\eq\label{thm1hpde}
h_t+\frac{1}{2}\,a\,h_{xx}=0,\quad h(T,\cdot)=v^{(-1)}
\en
in the classical sense, the superscript $(-1)$ denotes the spatial inverse and $B$ is a one-dimensional standard Brownian motion. 
\end{thm}

Our main result on the solvability of the system \eqref{cond_sys} requires the following definition.
\begin{defn}\label{cons_gen_def}
Let $M^\lambda$, $\lambda\in\Lambda$ be a family of Brownian martingales with time-dependent generators $\frac{1}{2}\,\sigma^\lambda(t,x)^2\frac{\mathrm{d}^2}{\mathrm{d}x^2}$ and initial values $m_0^\lambda$, respectively, and fix the notations $\Sigma^{\lambda}(t,x)=\int_{m_0^{\lambda}}^x \frac{\mathrm{d}y}{\sigma^\lambda(t,y)}$, $\lambda\in\Lambda$ and $\Theta^\lambda=(\Sigma^\lambda)^{(-1)}$, $\lambda\in\Lambda$. We call the family $M^\lambda$, $\lambda\in\Lambda$ \textit{consistent} if the following conditions hold. For any $\lambda_1\neq\lambda_2$ in $\Lambda$, there is a function $b^{\lambda_1,\lambda_2}:\,[0,T]\rightarrow\rr$, which solves the ordinary differential equation (ODE)
\eq\label{bODEmanyeq}
\begin{split}
\dot{b}^{\lambda_1,\lambda_2}=-\frac{1}{2}\sigma^{\lambda_1}_x(t,\Theta^{\lambda_1}(t,\Sigma^{\lambda_2}+b^{\lambda_1,\lambda_2}))
+\Sigma^{\lambda_1}_t(t,\Theta^{\lambda_1}(t,\Sigma^{\lambda_2}+b^{\lambda_1,\lambda_2}))\\
+\frac{1}{2}\sigma^{\lambda_2}_x-\Sigma^{\lambda_2}_t.
\end{split}
\en
for all values of $x$ such that the pointwise limit 
\eq\label{bTerminalmanyeq}
\Gamma^{\lambda_1,\lambda_2}:=\lim_{t\uparrow1} \Theta^{\lambda_1}(t,\Sigma^{\lambda_2}+b^{\lambda_1,\lambda_2})
\en
exists, belongs to ${\mathcal G}$, and one has
\eq\label{IVconsmanyeq}
m^{\lambda_2}_0=\ev^{m^{\lambda_1}_0}[(\Gamma^{\lambda_1,\lambda_2})^{(-1)}(M^{\lambda_1}(T))]. 
\en
Hereby, the argument of $\dot{b}^{\lambda_1,\lambda_2}$, $b^{\lambda_1,\lambda_2}$ is $t$, and the arguments of $\Sigma^{\lambda_2}$, $\sigma^{\lambda_2}_x$, $\Sigma^{\lambda_2}_t$ are $t$, $x$.  
\end{defn} 
It turns out that this notion of consistency is necessary and sufficient for the solvability of the system \eqref{cond_sys}. For a discussion of the conditions in Definition \ref{cons_gen_def}, please see Remarks \ref{cons_disc1} and \ref{cons_disc2} following the statement of Theorem \ref{thm2}.
 
\begin{thm}\label{thm2}
Suppose that $M^\lambda$, $\lambda\in\Lambda$ is a family of Brownian martingales, which all satisfy the conditions of Theorem \ref{thm1}. In addition, suppose that the diffusion coefficients $\sigma^\lambda$, $\lambda\in\Lambda$ of $M^\lambda$, $\lambda\in\Lambda$ are continuously differentiable and bounded away from $0$. Then, the system \eqref{cond_sys} is solvable by a Brownian martingale $X$ and a family $v^\lambda$, $\lambda\in\Lambda$ of functions in ${\mathcal G}$ if and only if the family $M^\lambda$, $\lambda\in\Lambda$ is consistent. Moreover, if this is the case, then there are uncountably infinitely many such solutions and the following are true:
\begin{enumerate}[(a)]
\item The processes $M^\lambda$, $\lambda\in\Lambda$ can be defined on the same probability space to form a weak solution of the degenerate system of SDEs
\eq\label{Mlambdasys}
\mathrm{d}M^\lambda(t)=\sigma^\lambda(t,M^\lambda(t))\,\mathrm{d}B(t),\quad M^\lambda(0)=m^\lambda_0
\en
on $[0,T]$. \\
\item If, in addition, there is a partial order $\prec$ on the set $\Lambda$ and the functions $\Gamma^{\lambda_1,\lambda_2}$ in Definition \ref{cons_gen_def} are such that $\lambda_1\prec\lambda_2$ implies $x\geq\Gamma^{\lambda_1,\lambda_2}(x)$ for all $x$ in the state space of $M^{\lambda_2}$, then the processes $M^\lambda$, $\lambda\in\Lambda$ can be defined on the same probability space in such a way that, for all $\lambda_1\prec\lambda_2$ in $\Lambda$, the inequality $M^{\lambda_1}(t)\leq M^{\lambda_2}(t)$ holds for all $t\in[0,T]$ with probability $1$. \\
\item If, in the situation of part (b), the set $\Lambda\subset\rr$ is an interval with the usual total order, the function $\lambda\mapsto m_0^\lambda$ is continuously differentiable and the function $\sigma(\lambda,t,x):=\sigma^\lambda(t,x)$ has continuous and bounded partial derivatives $\sigma_\lambda$, $\sigma_{\lambda\lambda}$, $\sigma_{\lambda\,x}$, $\sigma_x$ and $\sigma_{xx}$, then there is a version of the growth process $H(\lambda,t):=M^\lambda(t)$, $(\lambda,t)\in\Lambda\times[0,T]$, which is continuously differentiable in $\lambda$, and such that the corresponding derivative $H_\lambda$ solves 
\begin{eqnarray}\label{growth_dist}
&&\dot{H}_\lambda=\left(\sigma_\lambda(\lambda,t,H)+\sigma_x(\lambda,t,H)H_\lambda\right)\xi, \label{growth_eq} \\
&&H_\lambda(0,\lambda)=\frac{\mathrm{d}m^\lambda_0}{\mathrm{d}\lambda}.
\end{eqnarray}
Hereby, $\xi$ is a distribution-valued Gaussian field with correlation function 
\eq
\ev[\xi(t_1,\lambda_1)\xi(t_2,\lambda_2)]=\delta(t_2-t_1),\quad \lambda_1,\lambda_2\in\Lambda,\;t_1,t_2\in[0,T]. 
\en
\end{enumerate} 
\end{thm}

\smallskip

\begin{figure}[h] \label{fig1}
\begin{center}
\includegraphics[height=3.2in]{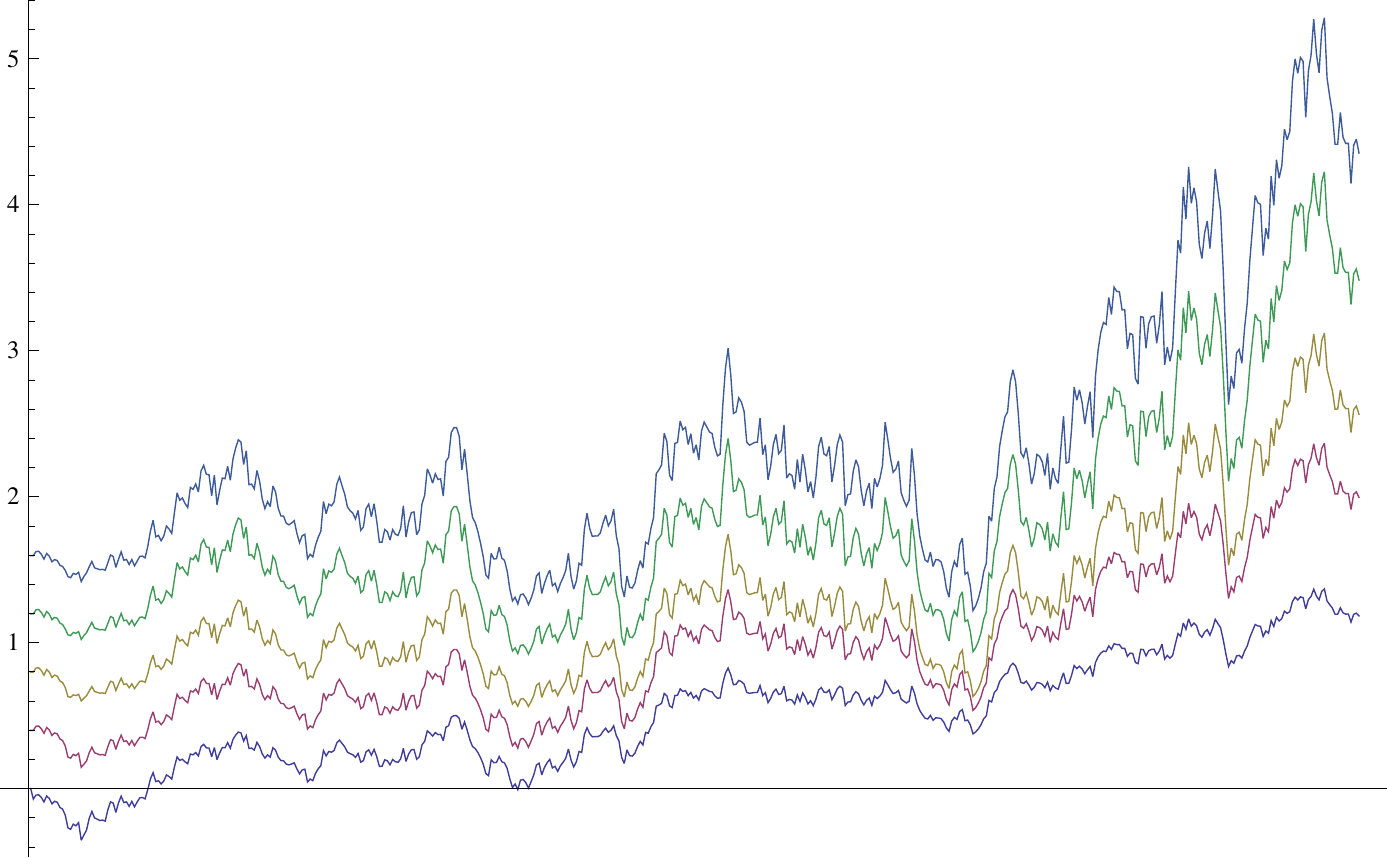}
\caption{The figure demonstrates the growth of a random surface as in Theorem \ref{thm2} (c). The family of consistent Brownian martingales is the one of Example \ref{ex_numeric} at the end of subsection \ref{subsec_systems} and the five curves correspond to sample paths of the growth process $H$ for $\lambda=1,2,3,4,5$.}
\label{capdist}
\end{center}
\end{figure}

\begin{rmk}\label{cons_disc1}
Consider first the situation $|\Lambda|=2$ and let $M$ be a Brownian martingale with diffusion coefficient $\sigma$ and initial value $m_0$, which satisfies the conditions of Theorem \ref{thm1}. Then, a family of Brownian martingales with diffusion coefficients $\tilde{\sigma}$ and initial values $\tilde{m}_0$, which are consistent with $M$, can be constructed as follows. First, let $\tilde{\Theta}$ be a solution of the backward heat equation
\begin{eqnarray*}
\tilde{\Theta}_t+\frac{1}{2}\,\tilde{\Theta}_{xx}
+\left(-\frac{1}{2}\,\sigma_x(t,\Sigma^{(-1)}(t,x+D))+\Sigma_t(t,\Sigma^{(-1)}(t,x+D))\right)\tilde{\Theta}_x=0 
\end{eqnarray*}
with some terminal condition $\tilde{\Theta}(T,\cdot)\in{\mathcal G}$, where $D$ is an arbitrary real constant. Next, set 
\eq
\tilde{\Sigma}=\tilde{\Theta}^{(-1)},\;\;\;\tilde{\sigma}=1/\tilde{\Sigma}_x\;\;\;\text{and}\;\;\;
\tilde{m}_0=\ev^{m_0}[\tilde{\Sigma}^{(-1)}(T,\Sigma(T,M(T))+D)]. 
\en
Then, the Brownian martingale $\tilde{M}$ corresponding to the diffusion coefficient $\tilde{\sigma}$ and initial value $\tilde{m}_0$ is consistent with $M$. Indeed, by differentiating the equation $\tilde{\Sigma}(t,\tilde{\Theta}(t,x))=x$ once with respect to $t$ and once and twice with respect to $x$, one can express the partial derivatives of $\tilde{\Theta}$ in terms of the partial derivatives of $\tilde{\Sigma}$. Plugging the resulting formulas into the backward heat equation and using $\tilde{\sigma}=1/\tilde{\Sigma}_x$, one then easily verifies that the constant function $b(t)=D$ solves the ODE in Definition \ref{cons_gen_def}. 

In the case that $\sigma\equiv1$, that is, when $M$ is a standard Brownian motion, the backward heat equation simplifies to
\eq
0=\tilde{\Theta}_t+\frac{1}{2}\,\tilde{\Theta}_{xx}.
\en
Next, take $\tilde{\Theta}(T,\cdot)$ to be an arbitrary function in ${\mathcal G}$. Then, $\tilde{\Theta}(T,\cdot)$ can be viewed as the cumulative distribution function of an infinite positive measure $\mu$ on $\rr$. Clearly, for any $t\in[0,T)$, the function $\tilde{\Theta}(t,\cdot)$ is given by a cumulative distribution function of the measure $\mu\ast\phi_{T-t}$, where $\phi_{T-t}$ is the normal density with mean $0$ and variance $T-t$. Consequently, $\tilde{\Sigma}(t,\cdot)=\tilde{\Theta}^{(-1)}(t,\cdot)$ is given by a quantile function $q_{T-t}$ of $\mu\ast\phi_{T-t}$. Thus, we conclude that the Brownian martingale solving
\eq
\mathrm{d}\tilde{M}(t)=\frac{1}{q'_{T-t}(\tilde{M}(t))}\,\mathrm{d}B(t)
\en  
is consistent with a Brownian motion started in $m_0$. We note hereby that, since we can choose the constant $D$ above in an arbitrary manner, we can take $\tilde{m}_0$ to be an arbitrary real number.   

Finally, in the case that $\Lambda$ has more than two elements, it suffices to observe that Theorems \ref{thm1} and \ref{thm2} imply that the definition of consistency for pairs of Brownian martingales induces a transitive relation on the space of Brownian martingales satisfying the conditions of Theorem \ref{thm1}. Indeed, if $M^1$, $M^2$, $M^3$ are Brownian martingales such that $M^1$ and $M^2$ are consistent and $M^2$ and $M^3$ are consistent, then by Theorem \ref{thm2}, there are functions $v^1,v^2,w^2,w^3\in{\mathcal G}$ and Brownian martingales $X$, $Y$ such that
\begin{eqnarray*}
&&M^1(\cdot)=\ev[v^1(X(T))|{\mathcal F}^X(\cdot)],\quad M^2(\cdot)=\ev[v^2(X(T))|{\mathcal F}^X(\cdot)],\\
&&M^2(\cdot)=\ev[w^2(Y(T))|{\mathcal F}^Y(\cdot)],\quad M^3(\cdot)=\ev[w^3(Y(T))|{\mathcal F}^Y(\cdot)].
\end{eqnarray*}  
Now, the uniqueness statement of Theorem \ref{thm1} implies that the processes $X$ and $(v^2)^{(-1)}(w^2(Y(\cdot)))$ coincide in law. Therefore, the process 
\[
\ev[v^1((v^2)^{(-1)}(w^2(Y(T))))|{\mathcal F}^Y(t)],\quad t\in[0,T] 
\]
has the law of $M^1$, so that $M^1$ and $M^3$ are consistent by Theorem \ref{thm2}. This shows the transitivity of the consistency property for pairs of Brownian martingales and, thus, that the family of processes $\tilde{M}$ constructed above is consistent, provided that they all satisfy the conditions of Theorem \ref{thm1}.   
\end{rmk}

\smallskip

\begin{rmk}\label{cons_disc2}
The following is the simplest example of a situation, in which the conditions of Definition \ref{cons_gen_def} hold. Suppose that for every $\lambda\in\Lambda$ the diffusion coefficient $\sigma^\lambda$ of $M^\lambda$ does not depend on $t$, and that for all $\lambda_1\neq\lambda_2$ in $\Lambda$ there are constants $A^{\lambda_1,\lambda_2}>0$, $B^{\lambda_1,\lambda_2}\in\rr$ such that
\[
\sigma^{\lambda_1}(x)=\frac{1}{A^{\lambda_1,\lambda_2}}\sigma^{\lambda_2}(A^{\lambda_1,\lambda_2}x+B^{\lambda_1,\lambda_2})
\]
and $m^{\lambda_2}_0=A^{\lambda_1,\lambda_2}m^{\lambda_1}_0+B^{\lambda_1,\lambda_2}$. Then, the family $M^\lambda$, $\lambda\in\Lambda$ is consistent. Indeed, choosing $b\equiv0$ and $\Gamma^{\lambda_1,\lambda_2}(x)=\frac{x-B^{\lambda_1,\lambda_2}}{A^{\lambda_1,\lambda_2}}$, one can check that the conditions of Definition \ref{cons_gen_def} are satisfied by using the identities $\Sigma^{\lambda_1}_{xx}=-\frac{\sigma^{\lambda_1}_x}{(\sigma^{\lambda_1})^2}$ and $\Sigma^{\lambda_2}_{xx}=-\frac{\sigma^{\lambda_2}_x}{(\sigma^{\lambda_2})^2}$.
\end{rmk}

\smallskip

\begin{rmk}
Considering the process $H$ in Theorem \ref{thm2} (c) as given, one can view the SPDE \eqref{growth_eq} for the partial derivative $H_\lambda$ as a degenerate version of the SPDEs studied in \cite{BM1}, \cite{BM2}, with random instead of deterministic coefficients. Indeed, setting ${\mathcal A}=0$ in equation (1.1) of \cite{BM1} and choosing the coefficient $g_1$ there to be the (random) diffusion coefficient in our equation \eqref{growth_eq} and the coefficient $f$ there to be the corresponding (random) It\^o correction term, one recovers our SPDE \eqref{growth_eq}. For other SPDEs with noise, which is white in time and correlated in space, we refer the reader to \cite{LS1}, \cite{LS2} and the references therein.       
\end{rmk}

\bigskip 
 
The rest of the paper is structured as follows. Section 2 is devoted to the proof of Theorem \ref{thm1}. As it turns out, equation \eqref{cond} can be solved by solving a nonlinear parabolic partial differential equation (PDE). This PDE is similar to the nonlinear PDEs appearing in \cite{mz} and can be reduced to a linear PDE by a suitable transformation. This is the content of subsection 2.1. In subsection 2.2, we solve the linear PDE and complete the proof of Theorem \ref{thm1}. In subsection 2.3, we analyze the solvability of equation \eqref{cond} in the degenerate case that $v=\mathbf{1}_{[c,\infty)}$ for some $c\in\rr$. In mathematical finance terms, the equation \eqref{cond} for this choice relates the price process $M$ of a digital option to its payoff $\mathbf{1}_{[c,\infty)}(X(T))$. Alternatively, one can view the process $M$ as describing the evolution of the conditional probabilities of the event $\{X(T)\geq c\}$ based on the information about the process $X$ available so far. Such processes can be observed in practice as quotes in prediction markets (see \cite{Al2} for more details). In contrast to the findings of Theorem \ref{thm1}, we find an uncountably infinite family of Brownian martingales solving the equation \eqref{cond} in this case. 

In section 3, we give the proof of Theorem \ref{thm2}. For the sake of clarity, we first prove the relevant statements of Theorem \ref{thm2} in the case that $|\Lambda|=2$. To this end, we need to solve a system of two nonlinear PDEs, each of the same type as in subsection 2.1. This is the content of subsection 3.1. In subsection 3.2, we show how the arguments extend to the case of general index sets $\Lambda$. Finally, in section 4 we give more examples of situations, in which Theorems \ref{thm1} and \ref{thm2} apply, by treating Brownian martingales consistent with the Kimura martingale.

\section{Proof of Theorem \ref{thm1}}\label{sec_1eq}

\subsection{Reduction to a linear PDE}

In this subsection, we show that the problem of solving \eqref{cond} naturally reduces to a backward Cauchy problem for a linear PDE. We start with a proposition.

\begin{prop}\label{prop3}
Let the functions $a$ and $v$ be as in Theorem \ref{thm1}, and suppose that $X$ is a Brownian martingale solving \eqref{cond} with this $v$ and $M$ being a Brownian martingale with time-dependent generator $\frac{1}{2}\,a(t,x)\,\frac{\mathrm{d}^2}{\mathrm{d}x^2}$. Moreover, assume that there exists a classical solution $u$ of 
\eq\label{u_pde}
u_t+\frac{1}{2}\frac{a(t,u)}{u_x^2}u_{xx}=0
\en
such that $u(T,\cdot)=v$, the partial derivative $u_{xxx}$ exists and is continuous, and
\eq\label{u_growth}
|u_x|\geq c_1 e^{-c_2|u|^2}\quad\text{for\;some\;constants\;}c_1,\,c_2>0.
\en 
Then, the time-dependent generator of $X$ is given by $\frac{1}{2}\,\frac{a(t,u(t,x))}{u_x(t,x)^2}\frac{\mathrm{d}^2}{\mathrm{d}x^2}$ and it holds
\eq
\ev[v(X(T))|{\mathcal F}^X(t)]=u(t,X(t)),\quad t\in[0,T]. 
\en
\end{prop} 

\smallskip

The proof of the proposition relies on the following lemma. 

\begin{lemma}\label{lemma4}
Let the functions $a$ and $v$ be as in Theorem \ref{thm1}, and suppose that $u$ is a classical solution of the PDE \eqref{u_pde} satisfying $u(T,\cdot)=v$ and \eqref{u_growth}, whose partial derivative $u_{xxx}$ exists and is continuous. Then, the spatial inverse $h:=u^{(-1)}$ is well-defined, solves the linear backward Cauchy problem 
\eq\label{h_pde}
h_t + \frac{1}{2}\,a(t,x)\,h_{xx}=0,\quad h(T,\cdot)=v^{(-1)}
\en 
in the classical sense and its partial derivative $h_x$ is bounded.
\end{lemma}

\medskip

\noindent\textbf{Proof.} From the PDE \eqref{u_pde} it is clear that either $u_x>0$ everywhere, or $u_x<0$ everywhere. However, since $v$ is strictly increasing, it holds $u_x(T,\cdot)=v'>0$ and, hence, the inequality $u_x>0$ must be true everywhere. Thus, for every $t\in[0,T]$, the function $u(t,\cdot)$ is strictly increasing and therefore the spatial inverse $h=u^{(-1)}$ is well-defined. Next, we consider the equation $u(t,h(t,x))=x$. Differentiating this equation with respect to $t$ once and with respect to $x$ once and twice, we obtain the following set of equations:
\begin{eqnarray} 
&&u_t(t,h(t,x))+u_x(t,h(t,x))h_t(t,x)=0, \label{uh1}\\
&&u_x(t,h(t,x))h_x(t,x)=1, \label{uh2}\\
&&u_{xx}(t,h(t,x))h_x(t,x)^2+u_x(t,h(t,x))h_{xx}(t,x)=0. \label{uh3}
\end{eqnarray}
Hereby, the partial derivatives $h_t$, $h_x$ and $h_{xx}$ exist by the Implicit Function Theorem. Plugging the equations \eqref{uh1}, \eqref{uh2}, \eqref{uh3} into \eqref{u_pde}, we obtain
\eq
-h_t(t,x)u_x(t,h(t,x))+\frac{1}{2}\frac{a(t,x)}{(1/h_x(t,x))^2}\frac{-u_x(t,h(t,x))h_{xx}(t,x)}{h_x(t,x)^2}=0.
\en
Simplifying, we conclude that $h$ is a classical solution of the problem \eqref{h_pde} as desired. Moreover, differentiating both sides of the equations in \eqref{h_pde} with respect to $x$, we deduce that $h_x$ is a classical solution of the problem 
\eq\label{hx_pde}
(h_x)_t+\frac{1}{2}\,a_x(t,x)\,(h_x)_x+\frac{1}{2}\,a(t,x)\,(h_x)_{xx}=0,\quad h_x(T,\cdot)=(v^{(-1)})_x. 
\en
Hereby, the existence and continuity of the partial derivative $h_{xxx}$ follows from the existence and continuity of $u_{xxx}$ and the Implicit Function Theorem. Moreover, the existence and continuity of the partial derivative $h_{tx}$ is a consequence of \eqref{h_pde} and the existence and continuity of $a_x$ and $h_{xxx}$. In addition, by assumption, the function $v'$ is bounded away from $0$, so that $h_x(T,\cdot)$ is bounded. Moreover, in view of \eqref{uh2}, the assumption \eqref{u_growth} implies 
\eq
|h_x|\leq \frac{1}{c_1}\,e^{c_2|x|^2}. 
\en
Applying the Maximum Principle for linear parabolic PDEs in the form of Theorem 9 in chapter 2 of \cite{fr}, we conclude that $h_x$ must be bounded everywhere. This finishes the proof of the lemma. \ep 

\bigskip

\noindent\textbf{Proof of Proposition \ref{prop3}.} By Lemma \ref{lemma4}, the spatial inverse $h=u^{(-1)}$ is well-defined, solves the PDE \eqref{h_pde} in the classical sense and has a bounded partial derivative $h_x$. Defining the process $N(t):=\ev[v(X(T))|{\mathcal F}^X(t)]$, $t\in[0,T]$, we obtain by It\^o's formula:
\eq\label{Nsde}
\mathrm{d}h(t,N(t))=h_x(t,N(t))\,\mathrm{d}N(t),\quad t\in(0,T). 
\en
Hereby, the drift terms disappear, since $N$ has the time-dependent generator $\frac{1}{2}\,a(t,x)\,\frac{\mathrm{d}^2}{\mathrm{d}x^2}$ by assumption and $h$ is a classical solution of the equation \eqref{h_pde}. Moreover, since $h_x$ is bounded and $N$ is a martingale, the process $h(t,N(t))$, $t\in[0,T]$ is also a martingale. Noting that $h(T,N(T))=u^{(-1)}(T,v(X(T)))=X(T)$ and recalling that $X$ is a martingale, we conclude that $h(t,N(t))=X(t)$ must hold for all $t\in[0,T]$ with probability $1$. Clearly, this implies that $N(t)=u(t,X(t))$, $t\in[0,T]$ with probability $1$. Hence, from equation \eqref{Nsde}, we see that the time-dependent generator of $X$ must be given by
\eq
\frac{1}{2}\,h_x(t,u(t,x))^2\,a(t,u(t,x))\,\frac{\mathrm{d}^2}{\mathrm{d}x^2}
=\frac{1}{2}\,\frac{a(t,u(t,x))}{u_x(t,x)^2}\,\frac{\mathrm{d}^2}{\mathrm{d}x^2}. 
\en    
Finally, by the definition of $N$ and the previous considerations, it holds
\eq
\ev[v(X(T))|{\mathcal F}^X(t)]=N(t)=u(t,X(t)),\quad t\in[0,T].
\en
This finishes the proof of the proposition. \ep

\bigskip

The proposition shows that, if a classical solution $u$ to the PDE \eqref{u_pde} with terminal condition $v$ exists and possesses a continuous partial derivative $u_{xxx}$, then a solution $X$ of the equation \eqref{cond} must have the time-dependent generator $\frac{1}{2}\,\frac{a(t,u(t,x))}{u_x(t,x)^2}\,\frac{\mathrm{d}^2}{\mathrm{d}x^2}$. In the next subsection, we show that, under the assumptions of Theorem \ref{thm1}, the desired classical solution of the nonlinear PDE \eqref{u_pde} exists by solving the \textit{linear} PDE \eqref{h_pde}, and that the resulting process $X$ is indeed a Brownian martingale. 

\subsection{Solution of the linear PDE}

As we have seen in the previous subsection, the crucial step in the solution of the equation \eqref{cond} consists of finding a sufficiently regular solution of the problem \eqref{h_pde}. In the following lemma we show that the latter exists under the assumptions of Theorem \ref{thm1}.

\begin{lemma}\label{lemma5}
Let the assumptions of Theorem \ref{thm1} hold. Then, the problem \eqref{h_pde} has a classical solution $h$, which satisfies
\eq\label{h_growth}
|h_x|\leq\frac{1}{c_1}e^{c_2|x|^2}\quad\text{for\;some\;constants\;}c_1,\,c_2>0
\en
and whose partial derivative $h_{xxx}$ exists and is continuous.
\end{lemma}

\medskip

\noindent\textbf{Proof.} Recalling the notation $p^M(t,x;T,y)$, $t\in[0,T)$, $x,y\in S$ for the transition densities of the process $M$, we define the function $h$ by 
\eq
h(t,x)=\int_S p^M(t,x;T,y)\,v^{(-1)}(y)\,\mathrm{d}y.
\en  
The assumptions on $p^M$ allow us to interchange the order of differentiation and integration to obtain
\begin{eqnarray}
&&h_t(t,x)=\int_S p^M_t(t,x;T,y)\,v^{(-1)}(y)\,\mathrm{d}y,\\
&&h_{xx}(t,x)=\int_S p^M_{xx}(t,x;T,y)\,v^{(-1)}(y)\,\mathrm{d}y,
\end{eqnarray}
and to conclude that the functions $h_t$ and $h_{xx}$ are continuous by using the Dominated Convergence Theorem. It follows that
\begin{eqnarray*}
&&h_t(t,x)+\frac{1}{2}\,a(t,x)\,h_{xx}(t,x)\\
&=&\int_S \big[p^M_t(t,x;T,y)+\frac{1}{2}\,a(t,x)\,p^M_{xx}(t,x;T,y)\big]\,v^{(-1)}(y)\,\mathrm{d}y=0
\end{eqnarray*}
is a consequence of the Kolmogorov backward equation
\eq
p^M_t(t,x;T,y)+\frac{1}{2}\,a(t,x)\,p^M_{xx}(t,x;T,y)=0.
\en
Moreover, by the definition of $h$, we have $h(T,\cdot)=v^{(-1)}$. In addition, the growth estimate \eqref{h_growth} follows from another exchange of the order of differentiation and integration 
\eq
h_x(t,x)=\int_S p^M_x(t,x;T,y)\,v^{(-1)}(y)\,\mathrm{d}y,
\en
$v^{(-1)}\in{\mathcal G}$ and \eqref{pMx_growth}. Finally, the partial derivative $h_{xxx}$ exists and is continuous due to yet another exchange of the order of differentiation and integration
\eq
h_{xxx}(t,x)=\int_S p^M_{xxx}(t,x;T,y)\,v^{(-1)}(y)\,\mathrm{d}y
\en
and the continuity of the right-hand side in the latter equation, justified by the assumptions on $p^M$ and the Dominated Convergence Theorem. \ep 

\bigskip

We are now ready to give the proof of Theorem \ref{thm1}.

\bigskip

\noindent\textbf{Proof of Theorem \ref{thm1}.} We first prove that a solution $X$ to \eqref{cond} of the described form exists. To this end, we let $h$ be a classical solution of the problem \eqref{h_pde} as in Lemma \ref{lemma5} and set $X(t)=h(t,M(t))$, $t\in[0,T]$. Applying It\^o's fomula, we see that
\eq\label{Xdyn}
\mathrm{d}X(t)=h_x(t,M(t))\,\mathrm{d}M(t),\quad t\in[0,T].
\en 
Hereby, the drift terms disappear, because $M$ has the time-dependent generator $\frac{1}{2}\,a(t,x)\,\frac{\mathrm{d}^2}{\mathrm{d}x^2}$ and $h$ is a classical solution of the PDE in \eqref{h_pde}. Moreover, noting that $h_x$ is a classical solution of the problem \eqref{hx_pde} and applying the Maximum Principle for linear parabolic equations in the form of Theorem 9 in chapter 2 of \cite{fr} (note the growth estimate \eqref{h_growth}), we conclude that $h_x$ is bounded between two positive constants. In particular, the spatial inverse $u:=h^{(-1)}$ is well-defined and $M(t)=u(t,X(t))$, $t\in[0,T]$. Putting these observations together with \eqref{Xdyn}, we conclude that $X$ is a Brownian martingale with time-dependent generator $\frac{1}{2}\,h_x(t,u(t,x))^2\,a(t,u(t,x))\,\frac{\mathrm{d}^2}{\mathrm{d}x^2}$. 

\bigskip

Next, we note that the identities $X(t)=h(t,M(t))$, $M(t)=u(t,X(t))$, $t\in[0,T]$ imply that the filtration generated by $M$ is the same as the filtration generated by $X$. Moreover, we have 
\eq
v(X(T))=v(v^{(-1)}(M(T)))=M(T). 
\en
Hence, since $M$ is a martingale in the filtration generated by $X$, it holds
\eq
M(t)=\ev[M(T)|{\mathcal F}^X(t)]=\ev[v(X(T))|{\mathcal F}^X(t)],\quad t\in[0,T]
\en
as desired.

\bigskip

We now turn to the proof of uniqueness. Differentiating both sides of the equation $h(t,u(t,x))=x$ once with respect to $t$ and twice with respect to $x$, we get 
\begin{eqnarray} 
&&h_t(t,u(t,x))+h_x(t,u(t,x))u_t(t,x)=0, \label{hu1}\\
&&h_{xx}(t,u(t,x))u_x(t,x)^2+h_x(t,u(t,x))u_{xx}(t,x)=0. \label{hu3}
\end{eqnarray}
Hence,
\[
u_t+\frac{1}{2}\frac{a(t,u)}{u_x^2}u_{xx}=-\frac{h_t(t,u)+\frac{1}{2}a(t,u)h_{xx}(t,u)}{h_x(t,u)}=0.
\]
In addition, the existence and continuity of the partial derivative $h_{xxx}$ and the Implicit Function Theorem imply that $u_{xxx}$ exists and is continuous. Moreover, the estimate \eqref{u_growth} is a direct consequence of the growth estimate \eqref{h_growth}. Applying Proposition \ref{prop3}, we conclude that a Brownian martingale $X$ solving \eqref{cond} must have the generator $\frac{1}{2}\,\frac{a(t,u(t,x))}{u_x(t,x)^2}\frac{\mathrm{d}^2}{\mathrm{d}x^2}$. This and the boundedness of the function $u_x$ show that the process $u(t,X(t))$, $t\in[0,T]$ is a martingale in the filtration generated by $X$ by an application of It\^o's formula. Therefore, by \eqref{cond} 
\eq\label{NX}
u(t,X(t))=\ev[v(X(T))|{\mathcal F}^X(t)]=N(t),\quad t\in[0,T]  
\en
for a process $N$ of the same law as $M$. From \eqref{NX} we see that the law of the process $X(t)=h(t,N(t))$, $t\in[0,T]$ is uniquely determined. \ep

\subsection{A degenerate case}

We turn now to the degenerate case that the function $v$ in equation \eqref{cond} is given by the indicator function $\mathbf{1}_{[c,\infty)}$ for some $c\in\rr$. To make sure that \eqref{cond} is well-posed, we need to assume that the state space of $M$ is the interval $[0,1]$ and that $M(T)\in\{0,1\}$ holds with probability $1$. As before, we write $\frac{1}{2}\,a(t,x)\,\frac{\mathrm{d}^2}{\mathrm{d}x^2}$ for the time-dependent generator of $M$. In contrast to the non-degenerate case of Theorem \ref{thm1}, there exist uncountably infinitely many Brownian martingales solving \eqref{cond} here.

\begin{prop}\label{prop_deg}
In the case that $v=\mathbf{1}_{[c,\infty)}$ for some $c\in\rr$ in equation \eqref{cond} and that $M$ is a Brownian martingale as just described, there exists an uncountably infinite family of solutions $X$ of equation \eqref{cond}, each of which satisfies an SDE of the form
\eq
\mathrm{d}X(t)=(c_1-c_0)\sqrt{a(t,(X(t)-c_0)/(c_1-c_0))}\,\mathrm{d}B(t),\quad t\in[0,T],
\en
for some constants $c_0<c\leq c_1$, where $B$ is a standard Brownian motion. 
\end{prop}

\medskip

\noindent\textbf{Proof.} As in the non-degenerate case, we consider the linear PDE
\eq\label{h_pde1}
h_t+\frac{1}{2}\,a(t,x)\,h_{xx}=0,
\en
which formally yields 
\[
h(t,x)=\ev[h(T,M(T))|M(t)=x]=x\,h(T,1)+(1-x)\,h(T,0). 
\]
This motivates letting $c_1:=h(T,1)\geq c$, $c_0:=h(T,0)<c$ and setting $X(t)=c_1 M(t)+c_0(1-M(t))$, $t\in[0,T]$. Indeed, then $X$ is a Brownian martingale. Moreover,
\[
\ev[\textbf{1}_{\{X(T)\geq c\}}|{\mathcal F}^X(t)]=\ev[M(T)|{\mathcal F}^X(t)]=M(t),
\]
since the filtration generated by $X$ is the same as the filtration generated by $M$, and $M$ is a martingale in this filtration. It follows that $X$ is a solution of the equation \eqref{cond}. Finally, by Theorem 4.2 in chapter 3 of \cite{KS}, the process $M$ obeys the SDE
\eq
\mathrm{d}M(t)=\sqrt{a(t,M(t))}\,\mathrm{d}B(t),\quad t\in[0,T]
\en
for a standard Brownian motion $B$. This shows that the process $X(t)=c_1 M(t)+c_0(1-M(t))$, $t\in[0,T]$ satisfies the SDE in the statement of the proposition. Now, it remains to note that the choice of the constants $c_0<c\leq c_1$ was arbitrary and the proof is finished. \ep
 
\section{Proof of Theorem \ref{thm2}}

\subsection{Systems of two equations}

For the sake of clarity of exposition and simpler notation, we first give the proof of the relevant statements of Theorem \ref{thm2} in the case that $|\Lambda|=2$. To fix notations, we let $M$, $\tilde{M}$ be two Brownian martingales with time-dependent generators ${\mathcal L}_t=\frac{1}{2}\,\sigma(t,m)^2\frac{\mathrm{d}^2}{\mathrm{d}m^2}$, $\tilde{{\mathcal L}}_t=\frac{1}{2}\,\tilde{\sigma}(t,\tilde{m})^2\frac{\mathrm{d}^2}{\mathrm{d}\tilde{m}^2}$ and initial values $m_0$, $\tilde{m}_0$, respectively. Assuming that $M$ and $\tilde{M}$ satisfy the conditions of Theorem \ref{thm1}, we seek Brownian martingales $X$ and functions $v,\tilde{v}\in{\mathcal G}$ such that the system of equations 
\begin{eqnarray}
&&M(\cdot)\overset{d}{=}\ev[v(X(T))|{\mathcal F}^X(\cdot)], \label{sys1}\\
&&\tilde{M}(\cdot)\overset{d}{=}\ev[\tilde{v}(X(T))|{\mathcal F}^X(\cdot)] \label{sys2} 
\end{eqnarray}
is satisfied. In the case that $|\Lambda|=2$, the definition of consistency (Definition \ref{cons_gen_def} in the introduction) simplifies to the following definition. 

\begin{defn}\label{cons2def}
Letting $\Sigma(t,x)=\int_{m_0}^x \frac{\mathrm{d}y}{\sigma(t,y)}$ and $\tilde{\Sigma}(t,x)=\int_{\tilde{m}_0}^x \frac{\mathrm{d}y}{\tilde{\sigma}(t,y)}$, we say that the Brownian martingales $M$ and $\tilde{M}$ are consistent if there is a function $b:\,[0,T]\rightarrow\rr$, which solves the ODE
\eq\label{bODEgen}
\dot{b}=-\frac{1}{2}\sigma_x(t,\Sigma^{(-1)}(t,\tilde{\Sigma}+b))+\Sigma_t(t,\Sigma^{(-1)}(t,\tilde{\Sigma}+b))
+\frac{1}{2}\tilde{\sigma}_x-\tilde{\Sigma}_t
\en
for all values of $x$ and such that the pointwise limit 
\eq\label{bTerminal}
\Gamma(\tilde{m}):=\lim_{t\uparrow T} \Sigma^{(-1)}(t,\tilde{\Sigma}+b)
\en
exists, belongs to ${\mathcal G}$, and it holds 
\eq\label{IVcons}
\tilde{m}_0=\ev^{m_0}[\Gamma^{(-1)}(M(T))]. 
\en
Hereby, the argument of $\dot{b}$, $b$ is $t$, and the arguments of $\tilde{\Sigma}$, $\tilde{\sigma}_x$, $\tilde{\Sigma}_t$ are $t$, $x$.
\end{defn}

\medskip

We show now that the system \eqref{sys1}, \eqref{sys2} is solvable if and only if the Brownian martingales $M$ and $\tilde{M}$ are consistent.
 
\begin{prop}\label{thm_sys2}
Suppose that the Brownian martingales $M$, $\tilde{M}$ satisfy the assumptions of Theorem \ref{thm1} and that their diffusion coefficients $\sigma$, $\tilde{\sigma}$ are continuously differentiable and bounded away from $0$. Then, there is a Brownian martingale $X$ and functions $v,\tilde{v}\in{\mathcal G}$, which solve the system \eqref{sys1}, \eqref{sys2}, if and only if $M$ and $\tilde{M}$ are consistent in the sense of Definition \ref{cons2def}. Moreover, if this is the case, then there are uncountably infinitely many such solutions.
\end{prop}

\medskip

\noindent\textbf{Proof.} We assume first that $M$ and $\tilde{M}$ are consistent in the sense of Definition \ref{cons2def}. In order to solve the system \eqref{sys1}, \eqref{sys2}, it suffices to find a Brownian martingale $X$ with a time-dependent generator $\frac{1}{2}\,a^X(t,x)\,\frac{\mathrm{d}^2}{\mathrm{d}x^2}$, as well as functions $v,\tilde{v}\in{\mathcal G}$ such that the system of PDEs
\begin{eqnarray}
&&u_t+\frac{1}{2}a^X\,u_{xx}=0, \label{uPDE}\\
&&\tilde{u}_t+\frac{1}{2}a^X\,\tilde{u}_{xx}=0, \label{vPDE}\\
&&u_x\sqrt{a^X}=\sigma(t,u), \label{uODE}\\ 
&&\tilde{u}_x\sqrt{a^X}=\tilde{\sigma}(t,\tilde{u}), \label{vODE}\\ 
&&u(T,\cdot)=v, \label{uIC}\\
&&\tilde{u}(T,\cdot)=\tilde{v} \label{vIC}
\end{eqnarray}
has a classical solution, $\ev[v(X(T))]=m_0$ and $\ev[\tilde{v}(X(T))]=\tilde{m}_0$. Indeed, then the processes $u(t,X(t))$, $t\in[0,T]$ and $\tilde{u}(t,X(t))$, $t\in[0,T]$ satisfy
\begin{eqnarray}
&&\mathrm{d}u(t,X(t))=\sigma(t,u(t,X(t)))\,\mathrm{d}B(t),\quad t\in[0,T], \label{uSDE}\\
&&\mathrm{d}\tilde{u}(t,X(t))=\tilde{\sigma}(t,\tilde{u}(t,X(t)))\,\mathrm{d}B(t),\quad t\in[0,T] \label{vSDE}
\end{eqnarray}
for a standard Brownian motion $B$ (apply It\^o's formula and \eqref{uPDE}, \eqref{uODE}, \eqref{vPDE}, \eqref{vODE}) and, thus, are Brownian martingales with time-dependent generators ${\mathcal L}_t=\frac{1}{2}\,\sigma(t,m)^2\frac{\mathrm{d}^2}{\mathrm{d}m^2}$, $\tilde{{\mathcal L}}_t=\frac{1}{2}\,\tilde{\sigma}(t,\tilde{m})^2\frac{\mathrm{d}^2}{\mathrm{d}\tilde{m}^2}$ and initial values $u(0,X(0))$, $\tilde{u}(0,X(0))$, respectively. Since 
\begin{eqnarray*}
&&\ev[v(X(T))|{\mathcal F}^X(\cdot)]=u(\cdot,X(\cdot)),\quad \ev[v(X(T))]=m_0,\\ 
&&\ev[\tilde{v}(X(T))|{\mathcal F}^X(\cdot)]=\tilde{u}(\cdot,X(\cdot)),\quad \ev[\tilde{v}(X(T))]=\tilde{m}_0,
\end{eqnarray*}
this gives the desired solution of the system \eqref{sys1}, \eqref{sys2}. Note hereby that by our assumptions on $a$ and $\tilde{a}$, the solutions to \eqref{uSDE} and \eqref{vSDE} are pathwise unique and, therefore, the martingale problems for ${\mathcal L}_t$ and $\tilde{{\mathcal L}}_t$ are well-posed due to the results of Yamada and Watanabe (see e.g. Proposition 3.20 in \cite{KS}).

\bigskip

To solve the system \eqref{uPDE}-\eqref{vIC}, we fix an arbitrary function $v\in{\mathcal G}$ and proceed as in the proof of Theorem \ref{thm1} to find a Brownian martingale $X$ and a function $u$ such that the equations \eqref{uPDE}, \eqref{uODE}, \eqref{uIC} and $m_0=\ev[v(X(T))]$ are satisfied. Next, to make sure that equation \eqref{vODE} holds, it suffices to choose $\tilde{u}$ in such a way that
\eq\label{cons1}
\frac{u_x}{\sigma(t,u)}=\frac{\tilde{u}_x}{\tilde{\sigma}(t,\tilde{u})}.
\en
However, recalling the definitions of the functions $\Sigma$ and $\tilde{\Sigma}$ in Definition \ref{cons2def} and integrating both sides of the equation \eqref{cons1} in $x$, we see that the latter is equivalent to the equation
\eq\label{cons2}
\Sigma(t,u)-\tilde{\Sigma}(t,\tilde{u})=b(t),\quad t\in[0,T]
\en
for some function $b:\,[0,T]\rightarrow\rr$. Moreover, since the functions $\sigma$ and $\tilde{\sigma}$ take only positive values by assumption, the functions $\Sigma(t,\cdot)$ and $\tilde{\Sigma}(t,\cdot)$ are strictly increasing for all $t\in[0,T]$. It is now easy to check that equation \eqref{cons2} will hold if and only if we set $\tilde{u}=\tilde{h}^{(-1)}$ with a function $\tilde{h}$ satisfying
\eq\label{whatistildeh}
\tilde{h}(t,x)=u^{(-1)}(t,\Sigma^{(-1)}(t,\tilde{\Sigma}(t,x)+b(t)))=:u^{(-1)}(t,r(t,x,b(t))).
\en 
Hereby, we recall from the proof of Theorem \ref{thm1} that the functions $u(t,\cdot)$, $t\in[0,T]$ are strictly increasing, so that the spatial inverses $u^{(-1)}(t,\cdot)$, $t\in[0,T]$ are well-defined and strictly increasing as well. We note at this point that the just introduced function $r$ is determined by the parameters of the problem. On the other hand, we are free to pick a function $b$ of our choice and choose it as the function $b$ in Definition \ref{cons2def}.

\bigskip

Next, we make sure that equation \eqref{vIC} holds. To this end, we note that by the definition of the function $\Gamma$ in Definition \ref{cons2def}
and \eqref{whatistildeh} we have
\eq
\tilde{h}(T,\cdot)=u^{(-1)}(T,\Gamma(\cdot))=v^{(-1)}(\Gamma(\cdot)).
\en
Thus, to ensure \eqref{vIC} it is enough to set $\tilde{v}=\Gamma^{(-1)}\circ v$. The resulting function $\tilde{v}$ belongs to ${\mathcal G}$ by our assumptions on $\Gamma$ and $v$. Moreover, one has
\eq
\ev[\tilde{v}(X(T))]=\ev[\Gamma^{(-1)}(v(X(T)))]=\ev^{m_0}[\Gamma^{(-1)}(M(T))]=\tilde{m}_0,
\en
where the second identity is a consequence of the fact that $X$, $v$ solve \eqref{sys1} by construction and the third identity is a consequence of $M$ and $\tilde{M}$ being consistent (see Definition \ref{cons2def}). 

\bigskip

It remains to check that with the choice of $\tilde{h}$ (and, hence, also of $\tilde{u}$) above, the equation \eqref{vPDE} is satisfied. We claim that it suffices to show that $\tilde{h}$ is a classical solution of the PDE
\eq\label{tildehpde}
\tilde{h}_t+\frac{1}{2}\,\tilde{\sigma}(t,x)^2\,\tilde{h}_{xx}=0.
\en
Indeed, differentiating the equation $\tilde{h}(t,\tilde{u}(t,x))=x$ once with respect to $t$ and once and twice with respect to $x$, one can express the partial derivatives of $\tilde{u}$ in terms of the partial derivatives of $\tilde{h}$ and compute
\[
\tilde{u}_t+\frac{1}{2}a^X\,\tilde{u}_{xx}=\tilde{u}_t+\frac{1}{2}\,\frac{\tilde{\sigma}(t,\tilde{u})^2}{\tilde{u}_x^2}\,\tilde{u}_{xx}
=-\frac{\tilde{h}_t(t,\tilde{u})}{\tilde{h}_x(t,\tilde{u})}-\frac{1}{2}\,\frac{\tilde{\sigma}(t,\tilde{u})^2 \tilde{h}_{xx}(t,\tilde{u})}{\tilde{h}_x(t,\tilde{u})}=0,
\]
where we have used \eqref{vODE} in the first identity. To check that $\tilde{h}$ is a classical solution of \eqref{tildehpde}, we set $h=u^{(-1)}$ and observe that, in view of \eqref{whatistildeh}, equation \eqref{tildehpde} can be rewritten as
\eq
\frac{\mathrm{d}}{\mathrm{d}t}h(t,r(t,x,b(t)))
+\frac{1}{2}\,\tilde{\sigma}(t,x)^2\frac{\mathrm{d}^2}{\mathrm{d}x^2}h(t,r(t,x,b(t)))=0.
\en
Simplifying the left-hand side, we see that this is equivalent to the following ODE for $b$:
\eq\label{bprelimODE}
\dot{b}(t)=-\frac{\frac{1}{2}\tilde{\sigma}(t,x)^2(h_{xx}r_x^2+h_x\,r_{xx})+h_t+h_x\,r_t}{r_b\,h_x}.
\en
Hereby, the arguments of $h$ and its partial derivatives are $t$, $r(t,x,b(t))$, and the arguments of $r$ and its partial derivatives are $t$, $x$, $b(t)$. The key observation is now that, due to the definitions of the functions $r$, $\Sigma$ and $\tilde{\Sigma}$ (see \eqref{whatistildeh} and Definition \ref{cons2def}), we have
\begin{eqnarray*}
\tilde{\sigma}(t,x)^2\,r_x(t,x,b)^2&=&\tilde{\sigma}(t,x)^2\,\Sigma^{(-1)}_x(t,\tilde{\Sigma}(t,x)+b)^2\,\tilde{\Sigma}_x(t,x)^2 \\
&=&\Sigma_x^{(-1)}(t,\tilde{\Sigma}(t,x)+b)^2=\sigma(t,r(t,x,b))^2.
\end{eqnarray*}
Plugging this into \eqref{bprelimODE} and recalling from the proof of Theorem \ref{thm1} that $h$ is a classical solution of the problem \eqref{h_pde}, we can simplify the ODE \eqref{bprelimODE} to
\eq\label{bprelimODE1} 
\dot{b}(t)=-\frac{\frac{1}{2}\tilde{\sigma}(t,x)^2\,r_{xx}(t,x,b(t))+r_t(t,x,b(t))}{r_b(t,x,b(t))}.
\en
Next, we evaluate the partial derivatives of $r$ to 
\begin{eqnarray*}
&&r_x=\frac{\tilde{\Sigma}_x}{\Sigma_x(t,\Sigma^{(-1)}(t,\tilde{\Sigma}+b))},\\
&&r_{xx}=\frac{\tilde{\Sigma}_{xx}}{\Sigma_x(t,\Sigma^{(-1)}(t,\tilde{\Sigma}+b))}
-\frac{\tilde{\Sigma}_x^2\,\Sigma_{xx}(t,\Sigma^{(-1)}(t,\tilde{\Sigma}+b))}{\Sigma_x(t,\Sigma^{(-1)}(t,\tilde{\Sigma}+b))^3}\\
&&r_t=\frac{-\Sigma_t(t,\Sigma^{(-1)}(t,\tilde{\Sigma}+b))+\tilde{\Sigma}_t}{\Sigma_x(t,\Sigma^{(-1)}(t,\tilde{\Sigma}+b))},\\
&&r_b=\frac{1}{\Sigma_x(t,\Sigma^{(-1)}(t,\tilde{\Sigma}+b))}.
\end{eqnarray*}
Putting this together with \eqref{bprelimODE1}, we end up with
\[
\dot{b}=-\frac{1}{2}\tilde{\sigma}^2\tilde{\Sigma}_{xx}
+\frac{1}{2}\tilde{\sigma}^2\frac{\Sigma_{xx}(t,\Sigma^{(-1)}(t,\tilde{\Sigma}+b))}{\Sigma_x(t,\Sigma^{(-1)}(t,\tilde{\Sigma}+b))^2}\tilde{\Sigma}_x^2
+\Sigma_t(t,\Sigma^{(-1)}(t,\tilde{\Sigma}+b))-\tilde{\Sigma}_t.
\]
Finally, recalling that $\Sigma_x=1/\sigma$ and $\tilde{\Sigma}_x=1/\tilde{\sigma}$, we can write the latter equation as
\begin{eqnarray*}
\dot{b}=\frac{1}{2}\sigma(t,\Sigma^{(-1)}(t,\tilde{\Sigma}+b))^2\,\Sigma_{xx}(t,\Sigma^{(-1)}(t,\tilde{\Sigma}+b))
+\Sigma_t(t,\Sigma^{(-1)}(t,\tilde{\Sigma}+b))\\
-\Big(\frac{1}{2}\tilde{\sigma}^2\,\tilde{\Sigma}_{xx}+\tilde{\Sigma}_t\Big).
\end{eqnarray*}
Moreover, since $\Sigma_{xx}=-\frac{\sigma_x}{\sigma^2}$ and $\tilde{\Sigma}_{xx}=-\frac{\tilde{\sigma}_x}{\tilde{\sigma}^2}$, this equation simplifies further to
\[
\dot{b}=-\frac{1}{2}\sigma_x(t,\Sigma^{(-1)}(t,\tilde{\Sigma}+b))+\Sigma_t(t,\Sigma^{(-1)}(t,\tilde{\Sigma}+b))
+\frac{1}{2}\tilde{\sigma}_x-\tilde{\Sigma}_t.
\]
The last equation holds due to the assumption that $M$ and $\tilde{M}$ are consistent. Thus, we have constructed a solution of the system \eqref{sys1}, \eqref{sys2}.

\bigskip

Conversely, suppose that $X$, $v$, $\tilde{v}$ form a solution of the system \eqref{sys1}, \eqref{sys2}. Then, by the uniqueness result in Theorem \ref{thm1}, the pair $(X,v)$ has to coincide with the solution of \eqref{sys1} constructed in the proof of Theorem \ref{thm1}, and the pair $(X,\tilde{v})$ has to coincide with the solution of \eqref{sys2} constructed in the proof of Theorem \ref{thm1}. In particular, it must hold
\eq
\frac{u_x}{\sigma(t,u)}=\frac{\tilde{u}_x}{\tilde{\sigma}(t,\tilde{u})}=\frac{1}{\sqrt{a^X}}.
\en 
One can now proceed as in the first part of the proof to deduce the existence of a function $b:\,[0,T]\rightarrow\rr$ such that $\tilde{h}$ is given by \eqref{whatistildeh}. Plugging this expression for $\tilde{h}$ into \eqref{tildehpde} and proceeding as before, one shows that $b$ must solve the ODE in Definition \ref{cons2def}. Moreover, setting $\Gamma=v\circ\tilde{v}^{(-1)}$, one easily verifies \eqref{bTerminal} and \eqref{IVcons} by using \eqref{cons2} and $\ev[\tilde{v}(X(T))]=\tilde{m}_0$, respectively. This shows the consistence of $M$ and $\tilde{M}$. 

\bigskip

Finally, since the choice of $v$ in the construction above was arbitrary among all functions in ${\mathcal G}$, we conclude that, if $M$ and $\tilde{M}$ are consistent, the system \eqref{sys1}, \eqref{sys2} has uncountably infinitely many solutions. \ep

\begin{rmk}\label{relax}
A careful reading of the proof of Proposition \ref{thm_sys2} shows that the system \eqref{sys1}, \eqref{sys2} is solvable by a Brownian martingale $X$ and functions $v,\tilde{v}\in{\mathcal G}$ under the weaker assumption that $M$ satisfies the conditions of Theorem \ref{thm1}, the martingale problem for $\frac{1}{2}\,\tilde{\sigma}(t,x)^2\,\frac{\mathrm{d}^2}{\mathrm{d}x^2}$ is well-posed and the diffusion coefficients $\sigma$, $\tilde{\sigma}$ are continuously differentiable, bounded and bounded away from $0$. The same is true for the statements in the upcoming Corollary \ref{thm_coupling}. 
\end{rmk}

\medskip

As a consequence of Proposition \ref{thm_sys2}, we obtain the existence of couplings of consistent Brownian martingales.

\begin{cor}\label{thm_coupling}
Suppose that the Brownian martingales $M$ and $\tilde{M}$ satisfy the conditions of Theorem \ref{thm1}, their diffusion coefficients are continuously differentiable and bounded away from $0$, and that $M$, $\tilde{M}$ are consistent in the sense of Definition \ref{cons2def}. Then, $M$ and $\tilde{M}$ can be defined on the same probability space to form a weak solution the degenerate system of SDEs
\begin{eqnarray}
&&\mathrm{d}M(t)=\sigma(t,M(t))\,\mathrm{d}B(t),\quad M(0)=m_0 \label{Msde}\\
&&\mathrm{d}\tilde{M}(t)=\tilde{\sigma}(t,\tilde{M}(t))\,\mathrm{d}B(t),\quad \tilde{M}(0)=\tilde{m}_0 \label{tildeMsde}
\end{eqnarray}
on $[0,T]$. If, in addition, the function $\Gamma$ in Definition \ref{cons2def} is such that $x\geq\Gamma(x)$ for all $x$ in the state space of $\tilde{M}$, then the inequality $M(t)\leq\tilde{M}(t)$ holds for all $t\in[0,T]$ with probability $1$. 
\end{cor}

\medskip

\noindent\textbf{Proof.} Applying Proposition \ref{thm_sys2}, we see that there is a probability space, a Brownian martingale $X$ defined on this space and functions $v,\tilde{v}\in{\mathcal G}$ such that the equations \eqref{sys1}, \eqref{sys2} hold. As we have seen in the course of the proof of Proposition \ref{thm_sys2}, the processes $\ev[v(X(T))|{\mathcal F}^X(\cdot)]$, $\ev[\tilde{v}(X(T))|{\mathcal F}^X(\cdot)]$ solve the degenerate system of SDEs \eqref{Msde}, \eqref{tildeMsde}. Moreover, under the additional assumption of $x\geq\Gamma(x)$ for all $x$ in the state space of $\tilde{M}$, we have
\eq
\ev[v(X(T))|{\mathcal F}^X(\cdot)]=\ev[\Gamma(\tilde{v}(X(T)))|{\mathcal F}^X(\cdot)]\leq\ev[\tilde{v}(X(T))|{\mathcal F}^X(\cdot)]
\en
with probability $1$, which is the desired monotone coupling. \ep

\begin{rmk}\label{strong}
In the setting of Remark \ref{cons_disc1} with $\sigma\equiv1$, Remark \ref{relax} and Corollary \ref{thm_coupling} show that the degenerate system of SDEs 
\begin{eqnarray}
&&\mathrm{d}M(t)=\mathrm{d}B(t), \label{BM_eq}\\
&&\mathrm{d}\tilde{M}(t)=\frac{1}{q'_{T-t}(\tilde{M}(t))}\,\mathrm{d}B(t) \label{quantile_eq}
\end{eqnarray}
has a weak solution on $[0,T]$ for any initial values $m_0$, $\tilde{m}_0$. We note hereby that the diffusion coefficient of $\tilde{M}$ is continuous, bounded and bounded away from $0$, so that the martingale problem satisfied by $\tilde{M}$ is well-posed. Moreover, since the weak solution constructed in the proof of Corollary \ref{thm_coupling} has the property that, for any $t\in[0,T]$, $\tilde{M}(t)$ can be written as a deterministic function of $M(t)$, the resulting process $\tilde{M}$ is adapted to the filtration generated by the Brownian motion $B$. This shows that the SDE \eqref{quantile_eq} has a strong solution for any initial value $\tilde{m}_0$. We note that if the function
\[
\Big(\frac{1}{q'_{T-t}(x)}\Big)'=(\tilde{\Theta}_x(t,q_{T-t}(x)))_x
=\frac{\tilde{\Theta}_{xx}(t,q_{T-t}(x))}{\tilde{\Theta}_x(t,q_{T-t}(x))}
=(\log\tilde{\Theta}_x)_x(t,q_{T-t}(x))
\]
fails to be bounded, the existence of a strong solution does not follow from classical existence theorems such as Theorem 2.9 in chapter 5 of \cite{KS}. In addition, in this case $\Gamma=q_0-D$, where $D$ is uniquely determined by
\[
\tilde{m}_0=\ev[\tilde{\Theta}(T,m_0+B(T)+D)].
\]
At this point, Corollary \ref{thm_coupling} shows that, whenever $q_0(x)\leq x+D$ holds for all $x\in\rr$, the system \eqref{BM_eq}, \eqref{quantile_eq} has a weak solution on $[0,T]$ starting from $(m_0,\tilde{m}_0)$ such that the inequality $m_0+B(t)\leq\tilde{M}(t)$ holds for all $t\in[0,T]$ with probability $1$. 
\end{rmk}

To demonstrate the possible range of applications of Corollary \ref{thm_coupling}, we give one immediate corollary.

\begin{cor}
Suppose that the Brownian martingales $M$ and $\tilde{M}$ satisfy the conditions of Theorem \ref{thm1}, their diffusion coefficients are continuously differentiable and bounded away from $0$, and that $M$, $\tilde{M}$ are consistent in the sense of Definition \ref{cons2def}. Assume further that the inequality $x\geq\Gamma(x)$ holds for all $x$ in the state space of $\tilde{M}$. Moreover, for any fixed $0<\epsilon<m_0$, let $\tau_\epsilon$, $\tilde{\tau}_\epsilon$ be the first hitting times of the set $[0,\epsilon]$ before $T$ by the respective Brownian martingales (which we set to be equal to $T$ if the set is not hit before $T$). Then, $\tau_\epsilon$ is stochastically dominated by $\tilde{\tau}_\epsilon$ in the sense that it holds 
\eq
\pp(\tau_\epsilon>t)\leq\pp(\tilde{\tau}_\epsilon>t)
\en
for all $t\in[0,T]$. 
\end{cor}

\subsection{Systems of any number of equations} \label{subsec_systems}

We now generalize the constructions of the previous subsection to give a proof of Theorem \ref{thm2} for a general index set $\Lambda$.

\bigskip

\noindent\textbf{Proof of Theorem \ref{thm2}.} We assume first that the family $M^\lambda$, $\lambda\in\Lambda$ is consistent in the sense of Definition \ref{cons_gen_def} and will construct a solution of the system \eqref{cond_sys} of the desired type. Arguing as in the proof of Proposition 
\ref{thm_sys2}, we see that it suffices to find a Brownian martingale $X$ with a time-dependent generator $\frac{1}{2}\,a^X(t,x)\,\frac{\mathrm{d}^2}{\mathrm{d}x^2}$, as well as functions $v^\lambda$, $\lambda\in\Lambda$ in ${\mathcal G}$ such that the system of PDEs
\begin{eqnarray}
&&u^\lambda_t+\frac{1}{2}a^X\,u^\lambda_{xx}=0,\quad\;\; \lambda\in\Lambda \label{uiPDE}\\
&&u^\lambda_x\sqrt{a^X}=\sigma^\lambda(t,u^\lambda),\quad \lambda\in\Lambda \label{uiODE}\\ 
&&u^\lambda(T,\cdot)=v^\lambda,\quad\quad\quad\;\; \lambda\in\Lambda \label{uiIC}
\end{eqnarray}
has a classical solution and $\ev[v^\lambda(X(T))]=m^\lambda_0$, $\lambda\in\Lambda$. To this end, we fix a $\lambda^*\in\Lambda$ and a function $v^{\lambda^*}\in{\mathcal G}$, and define $X$ as the solution of \eqref{cond} for the pair $(M^{\lambda^*},v^{\lambda^*})$, constructed in the proof of Theorem \ref{thm1}. Then, the equations \eqref{uiPDE}, \eqref{uiODE}, \eqref{uiIC} will hold for $\lambda^*$. Now, we set $v^\lambda=(\Gamma^{\lambda^*,\lambda})^{(-1)}\circ v^{\lambda^*}$ for all $\lambda\neq\lambda^*$ in $\Lambda$. Following the lines of the proof of Proposition \ref{thm_sys2}, one checks that with this choice the equations \eqref{uiPDE}, \eqref{uiODE}, \eqref{uiIC} are satisfied for all $\lambda\in\Lambda$ due to the consistency of the family $M^\lambda$, $\lambda\in\Lambda$. In addition, since there are uncountably infinitely many choices for the function $v^{\lambda^*}$, there are uncountably infinitely many solutions of the system \eqref{cond_sys}.

\bigskip

Conversely, suppose that a Brownian martingale $X$ and a family of functions $v^\lambda$, $\lambda\in\Lambda$ in ${\mathcal G}$ solve the system \eqref{cond_sys}. Then, by Proposition \ref{thm_sys2}, every pair $M^{\lambda_1}$, $M^{\lambda_2}$ must be consistent in the sense of Definition \ref{cons2def}. Therefore, the family $M^\lambda$, $\lambda\in\Lambda$ is consistent in the sense of Definition \ref{cons_gen_def}.  

\bigskip

At this point, to show the statements (a) and (b) in the theorem, one only needs to follow the lines of the proof of Corollary \ref{thm_coupling}. To prove statement (c), we apply Theorem 4.2 in \cite{Me} to deduce that, under the assumptions in statement (c) in the theorem, there exists a version of the unique strong solution of the system \eqref{Mlambdasys}, which is continuously differentiable in $\lambda$ for every $t\in[0,T]$, and, for all $\lambda\in\Lambda$, the derivative is a strong solution of the equation
\eq
\mathrm{d}\Big(\frac{\mathrm{d}M^\lambda}{\mathrm{d}\lambda}\Big)
=\left(\sigma_\lambda(\lambda,t,M^\lambda(t))+\sigma_x(\lambda,t,M^\lambda(t))\frac{\mathrm{d}M^\lambda}{\mathrm{d}\lambda}\right)\mathrm{d}B(t).
\en  
Rewriting the latter equation in distributional form, we arrive at \eqref{growth_dist}. \ep

\smallskip

\begin{rmk}\label{relax2}
A careful reading of the proof of Theorem \ref{thm2} shows that the system \eqref{cond_sys} has a solution of the desired type and the statements (a), (b), (c) of Theorem \ref{thm2} hold under the weaker assumption that there is a $\lambda^*\in\Lambda$ such that $M^{\lambda^*}$ satisfies the conditions of Theorem \ref{thm1}, the martingale problems for  $\frac{1}{2}\,\sigma^\lambda(t,x)^2\,\frac{\mathrm{d}^2}{\mathrm{d}x^2}$, $\lambda\in\Lambda\backslash\{\lambda^*\}$ are well-posed and the diffusion coefficients $\sigma^\lambda$, $\lambda\in\Lambda$ are continuously differentiable, bounded and bounded away from $0$.
\end{rmk}

\begin{rmk}\label{ex_growth}
Let $\Lambda\subset\rr$ be an interval and define the Brownian martingales $M^\lambda$, $\lambda\in\Lambda$ as in the second paragraph of Remark \ref{cons_disc1}, but choosing terminal conditions $\tilde{\Theta}^\lambda(T,\cdot)$ depending on $\lambda$. Letting $q^\lambda_{T-t}$ be the corresponding quantile functions, we recall from Remark \ref{cons_disc1} that the family $M^\lambda$, $\lambda\in\Lambda$ of Brownian martingales
with diffusion coefficients $\sigma^\lambda(t,x):=\frac{1}{(q^\lambda_{T-t})'(x)}$ is consistent for any choice of initial values $m^\lambda_0$, $\lambda\in\Lambda$. In view of Remark \ref{relax2} and Theorem \ref{thm2} (a), the degenerate system of SDEs 
\begin{eqnarray}
&&\mathrm{d}M(t)=\mathrm{d}B(t), \\
&&\mathrm{d}M^\lambda(t)=\frac{1}{(q^\lambda_{T-t})'(M^\lambda(t))}\,\mathrm{d}B(t),\quad \lambda\in\Lambda \label{quantile_eq_gen}
\end{eqnarray}
has a weak solution on $[0,T]$ for any initial values $m_0$ and $m^\lambda_0$, $\lambda\in\Lambda$. We recall hereby from Remark \ref{strong} that, for every $\lambda\in\Lambda$, the martingale problem satisfied by $M^\lambda$ is well-posed. Subsequently, we deduce as in Remark \ref{strong} that all processes $M^\lambda$, $\lambda\in\Lambda$ are adapted to the filtration generated by the Brownian motion $B$ by construction. This shows that the system \eqref{quantile_eq_gen} has a strong solution for any initial values $m^\lambda_0$, $\lambda\in\Lambda$. 

\smallskip

Next, define the constants $D^\lambda$, $\lambda\in\Lambda$ by 
\[
m^\lambda_0=\ev\big[\tilde{\Theta}^\lambda(T,m_0+B(T)+D^\lambda)\big]. 
\]
If the functions $\Gamma^\lambda:=q^\lambda_0-D^\lambda$ decrease pointwise in $\lambda$, then the functions $\Gamma^{\lambda_1,\lambda_2}:=(\Gamma^{\lambda_1})^{(-1)}\circ\Gamma^{\lambda_2}$ satisfy $\Gamma^{\lambda_1,\lambda_2}(x)\leq x$ for all $x$ in the state space of $M^{\lambda_2}$, whenever $\lambda_1\leq\lambda_2$. In this case, we can conclude from Theorem \ref{thm2} (b) that the system \eqref{quantile_eq_gen} has a weak solution, for which $M^{\lambda_1}(t)\leq M^{\lambda_2}(t)$ holds for all $t\in[0,T]$ with probability $1$, whenever $\lambda_1\leq\lambda_2$. Finally, if the functions $\lambda\mapsto m^\lambda_0$ and $(\lambda,t,x)\mapsto\frac{1}{(q_{T-t}^\lambda)'(x)}$ satisfy the additional regularity assumptions of Theorem \ref{thm2} (c), then there is a version of the growth process $H$ defined there, which satisfies 
\begin{eqnarray}\label{growth_dist_ex}
&&\dot{H}_\lambda=\left(\frac{-\big[\frac{(q^\lambda_{T-t})'}{\partial\lambda}(H)\big]}{(q^\lambda_{T-t})'(H)^2}
-\frac{(q^\lambda_{T-t})''(H)}{(q^\lambda_{T-t})'(H)^2}H_\lambda\right)\xi, \label{growth_eq_ex} \\
&&H_\lambda(0,\lambda)=\frac{\mathrm{d}m^\lambda_0}{\mathrm{d}\lambda}
\end{eqnarray}
with $\xi$ being the distribution-valued Gaussian field in part (c) of Theorem \ref{thm2}. 
\end{rmk}

\smallskip

\begin{exm}\label{ex_numeric}
Based on Remark \ref{ex_growth}, we now give a numeric example of a situation, in which Theorem \ref{thm2} (c) applies. To this end, we would like to choose the terminal conditions $\tilde{\Theta}^\lambda(T,\cdot)$ as the functions
\eq
f^\lambda(x):=x\,\mathbf{1}_{\{x\leq 0\}}+\lambda\,x\,\mathbf{1}_{\{x>0\}}
\en 
with $\lambda$ varying in $\Lambda:=[1,\infty)$. However, the functions $f^\lambda$ do not satisfy the differentiability assumption on $\tilde{\Theta}^\lambda(T,\cdot)$ of Remark \ref{ex_growth}. To avoid this problem, we choose a small smoothing parameter $\kappa\in(0,1)$ and set
\eq
\tilde{\Theta}^\lambda(T,\cdot):=f^\lambda\ast\phi_\kappa,
\en
where $\phi_\kappa$ is the normal density with mean $0$ and variance $\kappa$. A straightforward computation then gives
\begin{eqnarray*}
&&\tilde{\Theta}^\lambda(t,x)=\big(\tilde{\Theta}^\lambda(T,\cdot)\ast\phi_{T-t}\big)(x)=\big(f^\lambda\ast\phi_{T-t+\kappa}\big)(x)\\
&=&(\lambda-1)\sqrt{\frac{T-t+\kappa}{2\pi}}\exp\big(-x^2/(2(T-t+\kappa))\big)+x+(\lambda-1)x\,\Phi_{T-t+\kappa}(x),
\end{eqnarray*}
where $\phi_{T-t+\kappa}$ is the normal density with mean $0$ and variance $T-t+\kappa$ and $\Phi_{T-t+\kappa}$ is the corresponding cumulative distribution function. Therefore, the corresponding quantile functions $q^\lambda$, $\lambda\in[1,\infty)$ are given as the solutions of
\eq\label{qlambda_num}
\begin{split}
x=(\lambda-1)\sqrt{\frac{T-t+\kappa}{2\pi}}\,e^{-q^\lambda_{T-t}(x)^2/(2(T-t+\kappa))}+q^\lambda_{T-t}(x)\\
+(\lambda-1)q^\lambda_{T-t}(x)\,\Phi_{T-t+\kappa}(q^\lambda_{T-t}(x)).
\end{split}
\en
Differentiating both sides of this equation with respect to $x$, solving for $\frac{1}{(q^\lambda_{T-t})'}$ and simplifying, we end up with the identity
\eq
\frac{1}{(q^\lambda_{T-t})'(x)}=\frac{x-(\lambda-1)\sqrt{(T-t+\kappa)/(2\pi)}e^{-q^\lambda_{T-t}(x)^2/(2(T-t+\kappa))}}{q^\lambda_{T-t}(x)},
\en
where the functions $q^\lambda$, $\lambda\in[1,\infty)$ are given as the solutions of \eqref{qlambda_num}. We conclude from Remark \ref{ex_growth} that the Brownian martingales $M^\lambda$, $\lambda\in[1,\infty)$ corresponding to the diffusion coefficients $\sigma^\lambda(t,x)=\frac{1}{(q^\lambda_{T-t})'(x)}$, $\lambda\in[1,\infty)$ and initial values 
\begin{eqnarray*}
m^\lambda_0=\ev\big[\tilde{\Theta}^\lambda(T,B(T))\big]=\int_\rr\int_\rr f^\lambda(y)\phi_\kappa(x-y)\,\mathrm{d}y\,\phi_T(x)\,\mathrm{d}x \\
=\int_\rr f^\lambda(y)\phi_{T+\kappa}(y)\,\mathrm{d}y=(\lambda-1)\frac{\sqrt{T+\kappa}}{\sqrt{2\pi}},\quad\lambda\in[1,\infty) 
\end{eqnarray*}
can be defined on the same probability space in such a way that the inequalities $M^{\lambda_1}(t)\leq M^{\lambda_2}(t)$ are satisfied for all $t\in[0,T]$ and $\lambda_1\leq\lambda_2$ in $[1,\infty)$ with probability $1$. In Figure \ref{fig1}, we demonstrate this fact by plotting sample paths of the Brownian martingales $M^1,M^2,M^3,M^4,M^5$ for $\kappa=10^{-4}$ obtained by using the same sample path of the driving Brownian motion. 
\end{exm} 

\section{Additional examples}\label{examples_sec}

In this last section we demonstrate the results of Theorem \ref{thm1}, Proposition \ref{prop_deg} and Theorem \ref{thm2} on the example of the \textit{Kimura martingale} 
\eq\label{kimura}
\mathrm{d}M(t)=M(t)(1-M(t))\,\mathrm{d}B(t),\; t\in[0,T],\quad M(0)=m_0\in(0,1)
\en
and a time-changed version of it defined below. We recall that the state space of $M$ is the open interval $(0,1)$ and that $M$ does not exit from $(0,1)$ in finite time with probability $1$ (see e.g. section 6 in \cite{Hu}). 

\subsection{Single equation}

First, we fix a function $v\in{\mathcal G}$, and seek a Brownian martingale $X$ solving \eqref{cond} with $M$ given by \eqref{kimura}. To this end, we recall first that the transition densities of the Kimura martingale $M$ can be written down explicitly (see \cite{Hu}, \cite{Ki}):
\begin{eqnarray*}
p^M(t,x;T,y)=(2\pi(T-t))^{-1/2}\frac{(x(1-x))^{1/2}}{(y(1-y))^{3/2}}
\quad\quad\quad\quad\quad\quad\quad\quad\quad\quad\\
\times\exp\Big(-\frac{T-t}{8}-\frac{1}{2(T-t)}\Big(\log \frac{y(1-x)}{x(1-y)}\Big)^2\Big), 
\end{eqnarray*}
$t\in[0,T)$, $x,y\in(0,1)$. In particular, the partial derivatives $p_t^M$, $p_x^M$, $p^M_{xx}$ and $p^M_{xxx}$ exist and are continuous. In addition, a lengthy but straightforward calculation shows that
\eq
\lim_{y\downarrow0} p^M(t,x;T,y)=\lim_{y\uparrow1} p^M(t,x;T,y)=0
\en
and that the same is true for the partial derivatives $p^M_t$, $p^M_x$, $p^M_{xx}$ and $p^M_{xxx}$. Therefore, we can apply Theorem \ref{thm1} to conclude that there exists a unique Brownian martingale $X$, which satisfies \eqref{cond} for these $M$ and $v$. Moreover, Theorem \ref{thm1} and its proof imply that $X$ solves the SDE 
\eq
\mathrm{d}X(t)=h_x(t,h^{(-1)}(t,X(t)))\,h^{(-1)}(t,X(t))(1-h^{(-1)}(t,X(t)))\,\mathrm{d}B(t)
\en
on the time interval $[0,T]$, where 
\eq
h(t,x):=\int_{(0,1)} p^M(t,x;T,y)\,v^{(-1)}(y)\,\mathrm{d}y
\en
solves the problem \eqref{thm1hpde} with $a(t,x)=x^2(1-x)^2$ in the classical sense. 

\bigskip

To demonstrate an application of Proposition \ref{prop_deg}, we consider the process $\tilde{M}$ obtained from the Kimura martingale $M$ by the deterministic time change 
\eq
\alpha:\,[0,\infty)\rightarrow[0,T),\quad s\mapsto T(1-e^{-s}).
\en
In other words, $\tilde{M}$ is the Brownian martingale given by the solution of the SDE
\eq
\mathrm{d}\tilde{M}(t)=\frac{\tilde{M}(t)(1-\tilde{M}(t))}{\sqrt{T-t}}\,\mathrm{d}B(t),\;t\in[0,T),\quad \tilde{M}(0)=m_0\in(0,1). 
\en
Now, by Proposition 5.22 (d) in chapter 5 of \cite{KS}, the almost sure limit $\tilde{M}(T):=\lim_{t\uparrow T} \tilde{M}(t)=\lim_{s\rightarrow\infty} M(s)$ exists and is given by $1$ on a set of probability $m_0$ and by $0$ on a set of probability $1-m_0$. Letting $v=\mathbf{1}_{[c,\infty)}$ for some $c\in\rr$, we see from Proposition \ref{prop_deg} that the equation \eqref{cond} with $\tilde{M}$ and this choice of $v$ has uncountably infinitely many solutions, each of which satisfies an SDE of the form
\eq
\mathrm{d}X(t)=\frac{(X(t)-c_0)(c_1-X(t))}{(c_1-c_0)\sqrt{T-t}}\,\mathrm{d}B(t),\quad t\in[0,T],
\en
for some constants $c_0<c\leq c_1$.

\subsection{Systems of equations}

In order to give an example of an application of Theorem \ref{thm2}, we start with the Kimura martingale $M$ and proceed as in Remark \ref{cons_disc1} to construct a family of Brownian martingales consistent with $M$. In this particular case, one has
\begin{eqnarray}
&&\Sigma(t,x)=\log\frac{x(1-m_0)}{(1-x)m_0},\quad \Sigma^{(-1)}(t,x)=\frac{m_0\,e^x}{1-m_0+m_0\,e^x},\\
&&\Sigma_t(t,x)=0,\quad\sigma_x(t,x)=1-2x.
\end{eqnarray}
Hence, the backward heat equation of Remark \ref{cons_disc1} reads
\eq
0=\tilde{\Theta}_t+\frac{1}{2}\,\tilde{\Theta}_{xx}+\frac{1}{2}\,\frac{m_0\,e^{x+D}+m_0-1}{m_0\,e^{x+D}-m_0+1}\,\tilde{\Theta}_x.
\en
Setting $D=0$ and letting $\tilde{\Theta}(T,\cdot)$ be a function in ${\mathcal G}$, we see that, for any fixed $t\in[0,T]$, it holds $\tilde{\Theta}(t,\cdot)=\tilde{\Theta}(T,\cdot)\ast\psi_{T-t}$ with $\psi_{T-t}$ being the appropriate transition density of the diffusion
\eq\label{Rdiff}
\mathrm{d}R(t)=\frac{1}{2}\frac{m_0\,e^{R(t)}+m_0-1}{m_0\,e^{R(t)}-m_0+1}\,\mathrm{d}t+\mathrm{d}B(t). 
\en 
For each $t\in[0,T]$, we write $q_{T-t}$ for the spatial inverse of $\tilde{\Theta}(t,\cdot)$ and conclude from Remark \ref{cons_disc1} that the Kimura martingale $M$ is consistent with the Brownian martingale solving the SDE
\eq
\mathrm{d}\tilde{M}(t)=\frac{1}{q_{T-t}'(\tilde{M}(t))}\,\mathrm{d}B(t),
\quad \tilde{m}_0=\ev^{m_0}\Big[\tilde{\Theta}\Big(T,\log\frac{M(T)(1-m_0)}{(1-M(T))(m_0)}\Big)\Big]. 
\en

\smallskip

Letting $\Lambda\subset\rr$ be an interval, letting $\tilde{\Theta}(T,\cdot)$ depend on $\lambda$ and proceeding as before, we obtain a family of consistent Brownian martingales $M^\lambda$, $\lambda\in\Lambda$. Moreover, the family $M$, $M^\lambda$, $\lambda\in\Lambda$ satisfies the conditions described in Remark \ref{relax2}. Indeed, for every $\lambda\in\Lambda$, the boundedness, the boundedness away from $0$ and the continuity of the diffusion coefficients of $M^\lambda$ shows that the martingale problem satisfied by $M^\lambda$ is well-posed. Therefore, we may apply Theorem \ref{thm2} (a) in this situation to conclude that the degenerate system of SDEs
\begin{eqnarray*}\label{kimura_sys}
&&\mathrm{d}M(t)=M(t)(1-M(t))\,\mathrm{d}B(t),\\
&&\mathrm{d}M^\lambda(t)=\frac{1}{(q_{T-t}^\lambda)'(M^\lambda(t))}\,\mathrm{d}B(t),\;\;
m^\lambda_0=\ev^{m_0}\Big[\tilde{\Theta}^\lambda\Big(T,\log\frac{M(T)(1-m_0)}{(1-M(T))m_0}\Big)\Big]
\end{eqnarray*} 
has a weak solution on $[0,T]$. 

\bigskip

If, in addition, the functions $\tilde{\Theta}^\lambda(T,\cdot)$ are increasing pointwise in $\lambda$, then we can apply Theorem \ref{thm2} (b) to conclude that there is a weak solution of the latter system of SDEs on $[0,T]$ such that, for all $\lambda_1\leq\lambda_2$ in $\Lambda$, the inequality $M^{\lambda_1}(t)\leq M^{\lambda_2}(t)$ holds for all $t\in[0,T]$ with probability $1$. 

\bigskip

If, in addition to the above, the function $\lambda\mapsto m_0^\lambda$ is continuously differentiable and the function $\sigma(\lambda,t,x):=\frac{1}{(q^\lambda_{T-t})'(x)}$ has continuous and bounded partial derivatives $\sigma_\lambda$, $\sigma_{\lambda\lambda}$, $\sigma_{\lambda\,x}$, $\sigma_x$ and $\sigma_{xx}$, then by Theorem \ref{thm2} (c) there is a version of the growth process $H(\lambda,t)=M^\lambda(t)$, $\lambda\in\Lambda$, $t\in[0,T]$, which evolves according to
\begin{eqnarray}
&&\dot{H}_\lambda=\left(\frac{-\big[\frac{(q^\lambda_{T-t})'}{\partial\lambda}(H)\big]}{\big((q^\lambda_{T-t})'(H)\big)^2}
-\frac{(q^\lambda_{T-t})''(H)}{\big((q^\lambda_{T-t})'(H)\big)^2}H_\lambda\right)\xi, \\
&&H_\lambda(0,\lambda)=\frac{\mathrm{d}m^\lambda_0}{\mathrm{d}\lambda}
\end{eqnarray}
with $\xi$ being the distribution-valued Gaussian field in Theorem \ref{thm2}, part (c). Note that these equations coincide with the equations in Remark \ref{ex_growth} with the difference being that here the functions $q^\lambda_{T-t}$ are given by spatial inverses of convolutions with the heat kernel for the diffusion in \eqref{Rdiff}, as opposed to the heat kernel of a standard Brownian motion in Remark \ref{ex_growth}.  

\section{Acknowledgement}

The author would like to thank David J. Aldous for asking the question about the solvability of equation \eqref{cond}, which started this work, and also for his comments throughout the preparation of this work. He is also grateful to Amir Dembo for helpful suggestions. 

\bigskip

\bibliographystyle{alpha}

\begin{thebibliography}{50}

\bibitem{Al}
Aldous, D. (2012). Personal communication.

\bibitem{Al2}
Aldous, D. (2012). Using Prediction Market Data to Illustrate Undergraduate Probability. Preprint available at \textit{http://www.stat.berkeley.edu/~aldous/Papers/monthly.pdf}. 

\bibitem{BS}
Black, F., Scholes, M. (1973). The Pricing of Options and Corporate Liabilities. \textit{Journal of Political Economy} \textbf{81} 637-654.

\bibitem{BM1}
Buckdahn, R., Ma, J. (2001). Stochastic viscosity solutions for nonlinear stochastic partial differential equations. Part I. \textit{Stoch. Process Appl.} \textbf{93} 181-204.

\bibitem{BM2}
Buckdahn, R., Ma, J. (2001). Stochastic viscosity solutions for nonlinear stochastic partial differential equations. Part II. \textit{Stoch. Process Appl.} \textbf{93} 205-228. 

\bibitem{DR}
Dybvig, P., Ross, S. (1987). Arbitrage. \textit{Eatwell, J., Milgate, M., Newman, P. (eds.) The new Palgrave dictionary of economics} \textbf{1} 100-106. Macmillan, London. 

\bibitem{DS}
Delbaen, F., Schachermayer, W. (1994). A general version of the fundamental theorem of asset pricing. \textit{Math. Ann.} \textbf{300} 463-520.  

\bibitem{fr}
Friedman, A. (1964). \textit{Partial differential equations of parabolic type}. Prentice-Hall, Inc., Englewood Cliffs, N.J. 

\bibitem{HK}
Harrison, M., Kreps, D. (1979). Martingales and arbitrage in multiperiod security markets. \textit{J. Econ. Theory} \textbf{20} 381-408.

\bibitem{HP}
Harrison, M., Pliska, S. (1981). Martingales and stochastic integrals in the theory of continuous trading. \textit{Stoch. Process Appl.} \textbf{11} 215-260.

\bibitem{Hu} 
Huillet T. (2011). On the Karlin-Kimura approaches to the Wright-Fisher diffusion with fluctuating selection. \textit{J. Stat. Mech.} \textbf{P02016}.

\bibitem{KS}
Karatzas I., Shreve S. (1991). \textit{Brownian motion and stochastic calculus}. 2nd ed. Springer, New York.  

\bibitem{Ki} 
Kimura M. (1954). Process leading to quasi-fixation of genes in natural populations due to random fluctuations of selection intensities. \textit{Genetics} \textbf{39}. 

\bibitem{Kn}
Knight F. B. (1981). \textit{Essentials of Brownian motion and diffusion}. Mathematical Surveys and Monographs \textbf{18}. American Mathematical 
Society. 

\bibitem{LS1} 
Lions, P.-L., Souganidis, P. E. (1998). Fully nonlinear stochastic partial differential equations. \textit{C. R. Acad. Sci. Paris} \textbf{326} 1085-1092.

\bibitem{LS2}
Lions, P.-L., Souganidis, P. E. (1998). Fully nonlinear stochastic partial differential equations: non-smooth equations and applications. \textit{C. R. Acad. Sci. Paris} \textbf{327} 735-741.

\bibitem{Me}
Metivier M. (1983). Pathwise differentiability with respect to a parameter of solutions of stochastic differential equations. \textit{Lecture Notes in Control and Inform. Sci. Theory and application of random fields (Bangalore, 1982)} \textbf{49} 188-200.

\bibitem{Me2}
Merton, R. C. (1973). Theory of Rational Option Pricing. \textit{The Bell Journal of Economics and Management Science} \textbf{4} 141–-183.

\bibitem{mz}
Musiela M., Zariphopoulou, T. (2010). Stochastic partial differential equations and portfolio choice. \textit{Contemporary quantitative finance}  195-216. Springer, Berlin.
\end{thebibliography}

\bigskip

\end{document}